# Joint Modelling of Location, Scale and Skewness Parameters of the Skew Laplace Normal Distribution

Fatma Zehra Doğru[1*] and Olcay Arslan[2]

[1]Giresun University, Faculty of Economics and Administrative Sciences, Department of Econometrics, 28100 Giresun/Turkey
*fatma.dogru@giresun.edu.tr*
[2]Ankara University, Faculty of Science, Department of Statistics, 06100 Ankara/Turkey
*oarslan@ankara.edu.tr*

**Abstract**

In this article, we propose joint location, scale and skewness models of the skew Laplace normal (SLN) distribution as an alternative model for joint modelling location, scale and skewness models of the skew-t-normal (STN) distribution when the data set contains both asymmetric and heavy-tailed observations. We obtain the maximum likelihood (ML) estimators for the parameters of the joint location, scale and skewness models of the SLN distribution using the expectation-maximization (EM) algorithm. The performance of the proposed model is demonstrated by a simulation study and a real data example.

**Keywords:** EM algorithm, joint location, scale and skewness models, ML, SLN, SN.

## 1. Introduction

There are many remarkable and tractable methods for modeling the mean. In practice, modelling the dispersion will be of direct interest in its own right, to identify the sources of variability in the observations (Smyth and Verbyla (1999)).

In recent years, joint mean and dispersion models have been used for modeling heteroscedastic data sets. For instance, Park (1966) proposed a log linear model for the variance parameter and described the Gaussian model using a two stage process to estimate the parameters. Harvey (1976) examined the maximum likelihood (ML) estimation of the location and scale effects and also proposed a likelihood ratio test for heteroscedasticity. Aitkin (1987) proposed the modelling of variance heterogeneity in normal regression analysis. Verbyla (1993) estimated the parameters of the normal regression model under the log linear dependence of the variances on explanatory variables using the restricted ML. Engel and Huele (1996) represented an extension of the response surface approach to Taguchi type experiments for robust design by accommodating generalized linear modeling. Taylor and Verbyla (2004) introduced joint modelling of location and scale parameters of the t distribution. Lin and Wang (2009) proposed a robust approach for joint modelling of mean and scale parameters for longitudinal data. Lin and Wang (2011) studied Bayesian inference for joint modelling of location and scale parameters of the t distribution for longitudinal data. Wu and Li (2012) explored the variable selection for joint mean and dispersion models of the inverse Gaussian distribution. Li and Wu (2014) proposed joint modelling of location and scale parameters of the skew normal (SN) (Azzalini (1985, 1986)) distribution. Zhao and Zhang (2015) proposed variable selection of varying dispersion student-t regression models. Recently, Li et al. (2017) proposed variable selection in joint location, scale and skewness models of the SN distribution and Wu et al. (2017) explored variable selection in joint location, scale and skewness models of the STN distribution.

The skew exponential power distribution was proposed by Azzalini (1986) to deal with both skewness and heavy-tailedness, simultaneously. Its properties and inferential aspects were studied by DiCiccio and Monti (2004). Gómez et al. (2007) studied the skew Laplace normal (SLN) distribution that is a special case of the skew exponential power distribution. This distribution has wider range of skewness and also more applicable than the SN distribution. In literature, skewness and heavy-tailedness are modelled by using STN distribution for joint location, scale and skewness models. However, the STN distribution has an extra parameter that is the degrees of freedom parameter. Since this parameter should be estimated along with the other parameters, it may be computationally more exhaustive in practice. Therefore, in this paper, we propose to model joint location, scale and skewness models of the SLN distribution as an alternative model for the joint location, scale and skewness models of the STN distribution to model both skewness and heavy-tailedness in the data.

The rest of the paper is designed as follows. In Section 2, we give some properties of the SLN distribution. In Section 3, we introduce joint location, scale and skewness models of the SLN distribution. In Section 4, we give the ML estimation of the proposed joint location, scale and skewness model using the EM algorithm. In Section 5, we provide a simulation study to show the performance of the proposed model. In Section 6, modeling applicability of the proposed model is illustrated by using a real data set. The paper is finalized with a conclusion section.

## 2. Skew Laplace normal distribution

Let $Y$ be a SLN distributed random variable ($Y \sim SLN(\mu, \sigma^2, \lambda)$) with the location parameter $\mu \in \mathbb{R}$, scale parameter $\sigma^2 \in (0, \infty)$ and the skewness parameter $\lambda \in \mathbb{R}$. The probability density function (pdf) of $Y$ is given as

$$f(y) = 2f_L(y; \mu, \sigma)\Phi\left(\lambda \frac{y-\mu}{\sigma}\right), \tag{1}$$

where $f_L(y; \mu, \sigma)$ represents the pdf of Laplace distribution with

$$f_L(y; \mu, \sigma) = \frac{1}{2\sigma} e^{-\frac{|y-\mu|}{\sigma}}$$

and $\Phi$ is the cumulative distribution function of the standard normal distribution. Figure 1 displays the plots of the pdf of the SLN distribution for $\mu = 0$, $\sigma = 1$ and different values of $\lambda$.

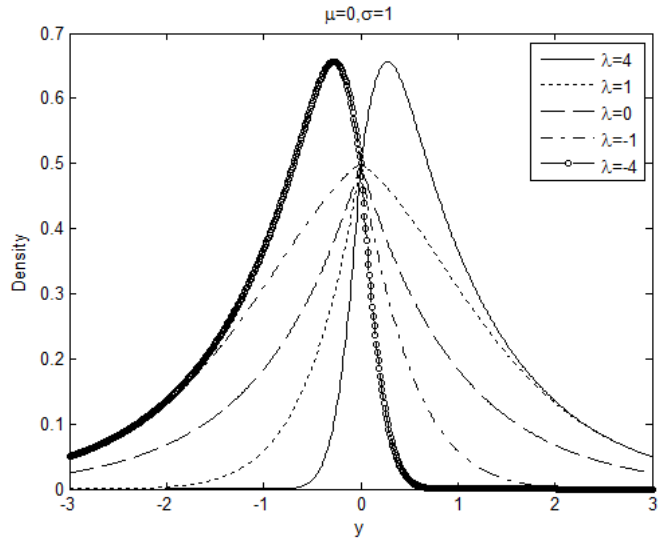

**Figure 1.** Examples of the SLN pdf for $\mu = 0$, $\sigma = 1$ and different skewness parameter values of $\lambda$.

Let the random variables $Z \sim SN(0,1,\lambda)$ and $V$ with the pdf $f_V(v) = v^{-3}\exp(-(2v^2)^{-1})$, $v > 0$ be two independent random variables. Then, the random variable $Y \sim SLN(\mu, \sigma^2, \lambda)$ has the following scale mixture form

$$Y = \mu + \sigma \frac{Z}{V}. \tag{2}$$

Further, using the stochastic representation of the SN (Azzalini (1986, p. 201) and Henze (1986, Theorem 1)) distributed random variable $Z$, the stochastic representation of the random variable $Y$ is obtained as

$$Y = \mu + \sigma \left( \frac{\lambda |Z_1|}{\sqrt{V^2(V^2 + \lambda^2)}} + \frac{Z_2}{\sqrt{V^2 + \lambda^2}} \right), \tag{3}$$

where $Z_1 \sim N(0,1)$ and $Z_2 \sim N(0,1)$ are independent random variables. This stochastic representation will give the following hierarchical representation of the SLN distribution. Let $U = \sqrt{V^{-2}(V^2 + \lambda^2)}|Z_1|$. Then,

$$\begin{aligned}
Y|u,v &\sim N\left(\mu + \frac{\sigma \lambda u}{v^2 + \lambda^2}, \frac{\sigma^2}{v^2 + \lambda^2}\right), \\
U|v &\sim TN\left(\left(0, \frac{v^2 + \lambda^2}{v^2}\right); (0, \infty)\right), \\
V &\sim f_V(v) = v^{-3} \exp(-(2v^2)^{-1}),
\end{aligned} \tag{4}$$

where $TN(\cdot)$ shows the truncated normal distribution. The hierarchical representation will allow us to carry on the parameter estimation using the EM algorithm. Using this hierarchical representation the joint pdf of $Y$, $U$ and $V$ can be written as

$$f(y,u,v) = \frac{1}{\pi \sigma} v^{-2} \exp(-(2v^2)^{-1}) \exp\left\{-\frac{1}{2}\left(\frac{v^2(y-\mu)^2}{\sigma^2} + \left(u - \frac{\lambda(y-\mu)}{\sigma}\right)^2\right)\right\}. \tag{5}$$

Next we will turn our attention to the conditional distribution of $U$ given $Y$ and $V$. Taking the integral of (5) over $U$, we obtain the joint pdf of $Y$ and $V$ as

$$f(y,v) = \left(\frac{2}{\pi \sigma^2}\right)^{1/2} v^{-2} \exp\left(-(2v^2)^{-1} - \frac{v^2 s^2}{2}\right) \Phi(\lambda s), \tag{6}$$

where $s = (y - \mu)/\sigma$. Then, dividing (5) by (6) yields the following conditional density function of $U$ given the others

$$f(u|y,v) = \frac{1}{\sqrt{2\pi}} \exp\left\{-\frac{(u - \lambda s)^2}{2}\right\} \Phi(\lambda s). \tag{7}$$

It is clear from the density given in (7) that $U$ and $V$ are conditionally independent. Therefore, the distribution of $U|Y = y$ is

$$U|Y = y \sim TN\big((\lambda s, 1); (0, \infty)\big). \tag{8}$$

Further, after dividing (6) by (1), we get the following conditional density function of $V$ given $Y$

$$f(v|y) = \sqrt{\frac{2}{\pi}} v^{-2} \exp\left\{-(2v^2)^{-1} - \frac{v^2 s^2}{2} + \frac{|y - \mu|}{\sigma}\right\}. \tag{9}$$

Now, we are ready to give the following proposition. The proof of this proposition can be easily done using the conditional pdfs given above.

**Proposition 1.** Using the hierarchical representation given in (4), we have the following conditional expectations

$$E(V^2|y) = \frac{\sigma}{|y-\mu|}, \qquad (10)$$

$$E(U|y) = \lambda s + \frac{\Phi(\lambda s)}{\phi(\lambda s)}, \qquad (11)$$

$$E(U^2|y) = 1 + \lambda s E(U|y). \qquad (12)$$

Note that these conditional expectations will be used in the EM algorithm given in Section 4.

### 3. Joint location, scale and skewness models of the SLN distribution

In this study, we consider the following joint location, scale and skewness models of the SLN distribution

$$\begin{cases} y_i \sim SLN(\mu_i, \sigma_i^2, \lambda_i), \; i = 1,2,\dots,n \\ \mu_i = x_i^T \boldsymbol{\beta}, \\ \log \sigma_i^2 = z_i^T \boldsymbol{\gamma}, \\ \lambda_i = w_i^T \boldsymbol{\alpha}, \end{cases} \qquad (13)$$

where $y_i$ is the $ith$ observed response, $\boldsymbol{x}_i = (x_{i1}, \dots, x_{ip})^T$, $\boldsymbol{z}_i = (z_{i1}, \dots, z_{iq})^T$ and $\boldsymbol{w}_i = (w_{i1}, \dots, w_{ir})^T$ are observed covariates corresponding to $y_i$, $\boldsymbol{\beta} = (\beta_1, \dots, \beta_p)^T$ is a $p \times 1$ vector of unknown parameters in the location model, and $\boldsymbol{\gamma} = (\gamma_1, \dots, \gamma_q)^T$ is a $q \times 1$ vector of unknown parameters in the scale model and $\boldsymbol{\alpha} = (\alpha_1, \dots, \alpha_r)^T$ is a $r \times 1$ vector of unknown parameters in the skewness model. These covariate vectors $\boldsymbol{x}_i, \boldsymbol{z}_i$ and $\boldsymbol{w}_i$ are not needed to be identical.

### 4. ML estimation of joint location, scale and skewness models of the SLN distribution

Let $(y_i, \boldsymbol{x}_i, \boldsymbol{z}_i), i = 1,2,\dots,n$, be a random sample from model given in (13). Let $\boldsymbol{\theta} = (\boldsymbol{\beta}, \boldsymbol{\gamma}, \boldsymbol{\alpha})$. Then, the log-likelihood function of $\boldsymbol{\theta}$ based on the observed data is written as

$$\ell(\boldsymbol{\theta}) = -\frac{1}{2} \sum_{i=1}^n z_i^T \boldsymbol{\gamma} - \sum_{i=1}^n \frac{|y_i - x_i^T \boldsymbol{\beta}|}{e^{z_i^T \boldsymbol{\gamma}/2}} + \sum_{i=1}^n \log \Phi(\kappa_i), \qquad (14)$$

where $\kappa_i = w_i^T \boldsymbol{\alpha} \frac{(y_i - x_i^T \boldsymbol{\beta})}{e^{z_i^T \boldsymbol{\gamma}/2}}$. The ML estimator of $\boldsymbol{\theta}$ can be found by maximizing the equation (14). We can see that the direct maximization of this function does not seem very tractable, so numerical algorithms may be needed to approximate the possible maximizer of this function. Since, the SLN distribution has a scale mixture form, the EM algorithm (Dempster et al. (1977)) can be implemented to obtain the ML estimator for $\boldsymbol{\theta}$. To simplify the steps of the EM algorithm, we will use the stochastic representation of the SLN distribution given in (3).

Let $V$ and $U$ be the latent variables. Using the hierarchical representation given in (4), or the model (13) we get the following hierarchical representation

$$Y_i|u_i, v_i \sim N\left(x_i^T\beta + \frac{e^{z_i^T\gamma/2}(w_i^T\alpha)u_i}{v_i^2 + (w_i^T\alpha)^2}, \frac{e^{z_i^T\gamma}}{v_i^2 + (w_i^T\alpha)^2}\right),$$

$$U_i|v_i = 1 \sim TN\left(\left(0, \frac{v_i^2 + (w_i^T\alpha)^2}{v_i^2}\right); (0, \infty)\right),$$

$$v_i \sim f(v_i) = v_i^{-3}\exp\left(-(2v_i^2)^{-1}\right). \tag{15}$$

Let $\boldsymbol{u} = (u_1, \ldots, u_n)$ and $\boldsymbol{v} = (v_1, \ldots, v_n)$ be the missing data and $(\boldsymbol{y}, \boldsymbol{u}, \boldsymbol{v})$ be the complete data, where $\boldsymbol{y} = (y_1, \ldots, y_n)$. Then, using the hierarchical representation given in (15), the complete data log-likelihood function of $\boldsymbol{\theta}$ can be written as

$$\ell_c(\boldsymbol{\theta}; \boldsymbol{y}, \boldsymbol{u}, \boldsymbol{v}) = \sum_{i=1}^{n}\left\{-\log\pi - \frac{1}{2}z_i^T\gamma - 2\log v_i - (2v_i^2)^{-1} - \frac{1}{2}\left(\frac{(y_i - x_i^T\beta)^2}{e^{z_i^T\gamma}}v_i^2\right.\right.$$

$$\left.\left. + u_i^2 - 2\frac{w_i^T\alpha}{e^{z_i^T\gamma/2}}(y_i - x_i^T\beta)u_i + \frac{(w_i^T\alpha)^2}{e^{z_i^T\gamma}}(y_i - x_i^T\beta)^2\right)\right\}. \tag{16}$$

To obtain the ML estimator of $\boldsymbol{\theta}$, we have to maximize (16). However, the estimators obtained from this maximization will be dependent on the latent variables. Thus, to handle this latency problem, we have to take the conditional expectation of the complete data log-likelihood function given the observed data $y_i$

$$E(\ell_c(\boldsymbol{\theta}; \boldsymbol{y}, \boldsymbol{u}, \boldsymbol{v})|y_i) = \sum_{i=1}^{n}\left\{-\log\pi - \frac{1}{2}z_i^T\gamma - 2E(\log V_i|y_i) - E\left((2V_i^2)^{-1}\big|y_i\right)\right.$$

$$-\frac{1}{2}\left(\frac{(y_i - x_i^T\beta)^2}{e^{z_i^T\gamma}}E(V_i^2|y_i) + E(U_i^2|y_i)\right.$$

$$\left.\left. -2\frac{w_i^T\alpha}{e^{z_i^T\gamma/2}}(y_i - x_i^T\beta)E(U_i|y_i) + \frac{(w_i^T\alpha)^2}{e^{z_i^T\gamma}}(y_i - x_i^T\beta)^2\right)\right\}. \tag{17}$$

The conditional expectations $E(V_i^2|y_i)$, $E(U_i|y_i)$ and $E(U_i^2|y_i)$ in (17) can be calculated using the conditional expectations given in (10)-(12). Note that since the other conditional expectations are not related to the parameters, we do not calculate them. Let

$$\hat{v}_i = E(V_i^2|y_i) = \frac{e^{z_i^T\hat{\gamma}/2}}{|y_i - x_i^T\hat{\beta}|}, \tag{18}$$

$$\hat{u}_{1i} = E(U_i|y_i) = \hat{\kappa}_i + \frac{\Phi(\hat{\kappa}_i)}{\phi(\hat{\kappa}_i)}, \tag{19}$$

$$\hat{u}_{2i} = E(U_i^2|y) = 1 + \hat{\kappa}_i\hat{u}_{1i}, \tag{20}$$

where, $\hat{\kappa}_i = w_i^T\hat{\alpha}\frac{(y_i - x_i^T\hat{\beta})}{e^{z_i^T\hat{\gamma}/2}}$. Then, using these conditional expectations in (17) we get the following objective function to be maximized with respect to $\boldsymbol{\theta}$

$$Q(\boldsymbol{\theta};\widehat{\boldsymbol{\theta}}) = \sum_{i=1}^{n}\left\{-\log\pi - \frac{1}{2}\mathbf{z}_i^T\boldsymbol{\gamma} - \frac{1}{2}\left(\frac{(y_i - \mathbf{x}_i^T\boldsymbol{\beta})^2}{e^{\mathbf{z}_i^T\boldsymbol{\gamma}}}\widehat{v}_i + \widehat{u}_{2i}\right.\right.$$
$$\left.\left. -2\frac{\mathbf{w}_i^T\boldsymbol{\alpha}}{e^{\mathbf{z}_i^T\boldsymbol{\gamma}/2}}(y_i - \mathbf{x}_i^T\boldsymbol{\beta})\widehat{u}_{1i} + \frac{(\mathbf{w}_i^T\boldsymbol{\alpha})^2}{e^{\mathbf{z}_i^T\boldsymbol{\gamma}}}(y_i - \mathbf{x}_i^T\boldsymbol{\beta})^2\right)\right\}. \quad (21)$$

Now, the steps of the EM algorithm will be as follows:

**EM algorithm:**

**1.** Take initial value for $\boldsymbol{\theta}^{(0)} = (\boldsymbol{\beta}^{(0)}, \boldsymbol{\gamma}^{(0)}, \boldsymbol{\alpha}^{(0)})$.

**2. E-Step:** Given the observed data and the current parameter values, find the conditional expectation of the complete data log-likelihood function given in (16). This corresponds to calculating the conditional expectations given in (18)-(20). This step will be carried on as follows. Compute the following conditional expectations for the $k = 0,1,2,...$ iteration

$$\widehat{v}_i^{(k)} = E(V_i^2|y_i, \widehat{\boldsymbol{\theta}}^{(k)}) = \frac{e^{\mathbf{z}_i^T\widehat{\boldsymbol{\gamma}}^{(k)}/2}}{|y_i - \mathbf{x}_i^T\widehat{\boldsymbol{\beta}}^{(k)}|}, \quad (22)$$

$$\widehat{u}_{1i}^{(k)} = E(U_i|y_i, \widehat{\boldsymbol{\theta}}^{(k)}) = \widehat{\kappa}_i^{(k)} + \frac{\Phi(\widehat{\kappa}_i^{(k)})}{\phi(\widehat{\kappa}_i^{(k)})}, \quad (23)$$

$$\widehat{u}_{2i}^{(k)} = E(U^2|y, \widehat{\boldsymbol{\theta}}^{(k)}) = 1 + \widehat{\kappa}_i^{(k)}\widehat{u}_{1i}^{(k)}, \quad (24)$$

where, $\widehat{\kappa}_i^{(k)} = \mathbf{w}_i^T\widehat{\boldsymbol{\alpha}}^{(k)}\frac{(y_i - \mathbf{x}_i^T\widehat{\boldsymbol{\beta}}^{(k)})}{e^{\mathbf{z}_i^T\widehat{\boldsymbol{\gamma}}^{(k)}/2}}$.

**3. M-Step:** Use these conditional expectations in $Q(\boldsymbol{\theta};\widehat{\boldsymbol{\theta}})$ and maximize it with respect to $\boldsymbol{\theta}$ to obtain new estimates. This maximization step yields the following formulation to update the new estimates. The $(k+1)th$ parameter estimates can be computed using

$$\widehat{\boldsymbol{\theta}}^{(k+1)} = \widehat{\boldsymbol{\theta}}^{(k)} + \left(-H(\widehat{\boldsymbol{\theta}}^{(k)})\right)^{-1}G(\widehat{\boldsymbol{\theta}}^{(k)}), \quad (25)$$

where $G(\boldsymbol{\theta}) = \frac{\partial Q(\boldsymbol{\theta};\widehat{\boldsymbol{\theta}})}{\partial \boldsymbol{\theta}}$ and $H(\boldsymbol{\theta}) = \frac{\partial^2 Q(\boldsymbol{\theta};\widehat{\boldsymbol{\theta}})}{\partial \boldsymbol{\theta} \partial \boldsymbol{\theta}^T}$.

**4.** Repeat E and M steps until the convergence is satisfied.

**Remark.** For the detail expressions of $G(\boldsymbol{\theta})$ and $H(\boldsymbol{\theta})$ see Appendix.

## 5. Simulation study

In this section, we give a simulation study to show the performance of the proposed location, scale and skewness models of the SLN distribution in terms of mean squared error (MSE). The MSE is given with the following formula

$$\widehat{MSE}(\widehat{\theta}) = \frac{1}{N}\sum_{j=1}^{N}(\widehat{\theta}_j - \theta)^2,$$

where $\theta$ is the true parameter value, $\hat{\theta}_j$ is the estimate of $\theta$ for the $jth$ simulated data and $\bar{\theta} = \frac{1}{N}\sum_{j=1}^{N}\hat{\theta}_j$. All simulation studies are conducted as $N = 1000$ times. We set the sample sizes as $50, 100, 150$ and $200$. Note that the simulation study and real data example are performed using MATLAB R2015b. For all numerical calculations, the convergence rule is taken as $10^{-6}$.

The data are generated from the following location, scale and skewness models of the SLN distribution

$$\begin{cases} y_i \sim SLN(\mu_i, \sigma_i^2, \lambda), & i = 1,2,\ldots,n \\ \mu_i = x_i^T \beta, \\ \log \sigma_i^2 = z_i^T \gamma, \\ \lambda_i = w_i^T \alpha. \end{cases}$$

Here, all covariate vectors $x_i, z_i$ and $w_i$ are independently generated from uniform distribution $U(-1,1)$. To carry out the simulation study, we take the following two cases for true parameter values:

Case I: $\beta_0 = (0,-1,-1)^T, \gamma_0 = (0,-1,-1)^T$ and $\alpha_0 = (0,-1,-1)^T$,
Case II: $\beta_0 = (0,1,1)^T, \gamma_0 = (0,1,1)^T$ and $\alpha_0 = (0,1,1)^T$,
Case III: $\beta_0 = (1,1,0,0,1)^T, \gamma_0 = (0.7,0.7,0,0,0.7)^T$ and $\alpha_0 = (0.5,0.5,0,0,0.5)^T$.

The simulation results are summarized in Tables 1, 2 and 3. These tables include the mean of the estimators and the values of MSE. From these tables, we observe the followings. The proposed EM algorithm is working accurately for estimating the parameters. When the sample sizes increase, the values of MSE decrease.

**Table 1.** Mean of the estimators and the values of MSE for the different sample sizes for the Case I.

| Model | $n$ | 50 | | 100 | | 150 | | 200 | |
|---|---|---|---|---|---|---|---|---|---|
| | | Mean | MSE | Mean | MSE | Mean | MSE | Mean | MSE |
| Location Model | $\beta_0$ | 0.0010 | 0.0433 | 0.0020 | 0.0161 | -0.0002 | 0.0088 | -0.0006 | 0.0061 |
| | $\beta_1$ | -1.0068 | 0.0624 | -0.9983 | 0.0243 | -0.9987 | 0.0143 | -0.9965 | 0.0109 |
| | $\beta_2$ | -0.9950 | 0.0643 | -0.9987 | 0.0256 | -1.0077 | 0.0150 | -0.9987 | 0.0117 |
| Scale Model | $\gamma_0$ | -0.0655 | 0.0905 | -0.0384 | 0.0369 | -0.0278 | 0.0228 | -0.0223 | 0.0161 |
| | $\gamma_1$ | -1.0595 | 0.3300 | -1.0576 | 0.1376 | -1.0264 | 0.0829 | -1.0088 | 0.0594 |
| | $\gamma_2$ | -1.1024 | 0.3445 | -1.0281 | 0.1272 | -1.0118 | 0.0843 | -1.0060 | 0.0629 |
| Skewness Model | $\alpha_0$ | -0.0527 | 0.5398 | 0.0035 | 0.0600 | -0.0034 | 0.0296 | -0.0016 | 0.0190 |
| | $\alpha_1$ | -1.5780 | 1.7814 | -1.1981 | 0.2234 | -1.1213 | 0.1113 | -1.0852 | 0.0727 |
| | $\alpha_2$ | -1.6140 | 2.0557 | -1.2066 | 0.2339 | -1.1157 | 0.1141 | -1.0754 | 0.0652 |

**Table 2.** Mean of the estimators and the values of MSE for the different sample sizes for the Case II.

| Model | $n$ | 50 | | 100 | | 150 | | 200 | |
|---|---|---|---|---|---|---|---|---|---|
| | | Mean | MSE | Mean | MSE | Mean | MSE | Mean | MSE |
| Location Model | $\beta_0$ | 0.0059 | 0.0417 | 0.0035 | 0.0142 | -0.0020 | 0.0088 | 0.0018 | 0.0067 |
| | $\beta_1$ | 1.0046 | 0.0623 | 1.0040 | 0.0247 | 0.9956 | 0.0158 | 1.0034 | 0.0115 |
| | $\beta_2$ | 1.0036 | 0.0661 | 1.0095 | 0.0265 | 1.0085 | 0.0138 | 1.0036 | 0.0105 |
| Scale Model | $\gamma_0$ | -0.0794 | 0.0935 | -0.0369 | 0.0346 | -0.0209 | 0.0224 | -0.0239 | 0.0162 |
| | $\gamma_1$ | 1.0957 | 0.3421 | 1.0380 | 0.1282 | 1.0281 | 0.0796 | 1.0200 | 0.0600 |
| | $\gamma_2$ | 1.0562 | 0.3332 | 1.0357 | 0.1252 | 1.0193 | 0.0783 | 1.0130 | 0.0580 |
| Skewness Model | $\alpha_0$ | -0.0130 | 0.4643 | 0.0045 | 0.0635 | -0.0009 | 0.0285 | -0.0013 | 0.0203 |
| | $\alpha_1$ | 1.5993 | 1.9956 | 1.2105 | 0.2636 | 1.1210 | 0.1035 | 1.0717 | 0.0672 |
| | $\alpha_2$ | 1.5928 | 1.7924 | 1.2117 | 0.2621 | 1.1148 | 0.1075 | 1.0802 | 0.0705 |

**Table 3.** Mean of the estimators and the values of MSE for the different sample sizes for the Case III.

| Model | $n$ | 50 | | 100 | | 150 | | 200 | |
|---|---|---|---|---|---|---|---|---|---|
| | | Mean | MSE | Mean | MSE | Mean | MSE | Mean | MSE |
| Location Model | $\beta_0$ | 1.0035 | 0.1734 | 0.9631 | 0.0900 | 0.9596 | 0.0508 | 0.9808 | 0.0356 |
| | $\beta_1$ | 1.0018 | 0.2093 | 1.0057 | 0.0709 | 0.9966 | 0.0415 | 1.0012 | 0.0296 |
| | $\beta_2$ | 0.0117 | 0.2129 | -0.0181 | 0.0806 | 0.0045 | 0.0473 | 0.0091 | 0.0294 |
| | $\beta_3$ | -0.0151 | 0.2219 | -0.0006 | 0.0784 | 0.0037 | 0.0447 | -0.0095 | 0.0301 |
| | $\beta_4$ | 0.9971 | 0.2102 | 1.0053 | 0.0767 | 0.9958 | 0.0372 | 1.0081 | 0.0287 |
| Scale Model | $\gamma_0$ | 0.5812 | 0.1404 | 0.6667 | 0.0598 | 0.6913 | 0.0387 | 0.6792 | 0.0240 |
| | $\gamma_1$ | 0.8583 | 0.5997 | 0.7411 | 0.1412 | 0.7145 | 0.0955 | 0.7226 | 0.0685 |
| | $\gamma_2$ | -0.0148 | 0.6155 | -0.0024 | 0.1715 | 0.0170 | 0.0963 | -0.0147 | 0.0612 |
| | $\gamma_3$ | -0.0303 | 0.5936 | -0.0136 | 0.1621 | -0.0168 | 0.0951 | -0.0110 | 0.0665 |
| | $\gamma_4$ | 0.8151 | 0.5718 | 0.7538 | 0.1561 | 0.7044 | 0.0920 | 0.7309 | 0.0621 |
| Skewness Model | $\alpha_0$ | 1.3635 | 3.2188 | 0.7983 | 0.4941 | 0.6819 | 0.1866 | 0.6067 | 0.0975 |
| | $\alpha_1$ | 1.2163 | 3.1849 | 0.7769 | 0.4341 | 0.6210 | 0.1381 | 0.5894 | 0.0831 |
| | $\alpha_2$ | 0.0499 | 1.8253 | 0.0156 | 0.1861 | -0.0059 | 0.0807 | 0.0033 | 0.0420 |
| | $\alpha_3$ | -0.0369 | 1.6726 | -0.0245 | 0.2151 | -0.0103 | 0.0764 | -0.0041 | 0.0499 |
| | $\alpha_4$ | 1.3361 | 3.4589 | 0.7471 | 0.3766 | 0.6317 | 0.1347 | 0.5827 | 0.0608 |

## 6. Real data example

The Martin Marietta data set includes the relationship of the excess rate of returns of the Marietta Company and an index for the excess rate of return for the New York Exchange (CRSP). These rate of returns for the company and the CRSP index were determined monthly over a period of five years. This data set used by Butler et al. (1990) for modelling a simple linear regression with Gaussian errors. Azzalini and Capitanio (2003) analyzed this data set for modelling the linear regression model when the errors have the skew t distribution. Also, Taylor and Verbyla (2004) examined this data set for joint modelling of location and scale parameters of the t distribution. We display the scatter plot of the data set and the histogram of the Martin Marietta excess returns. Since the skewness coefficient of Martin Marietta excess returns is 2.9537 and also according to the Figure 2 (b), we can say that it will be more suitable to model this data set with a joint location, scale and skewness models of a skew distribution.

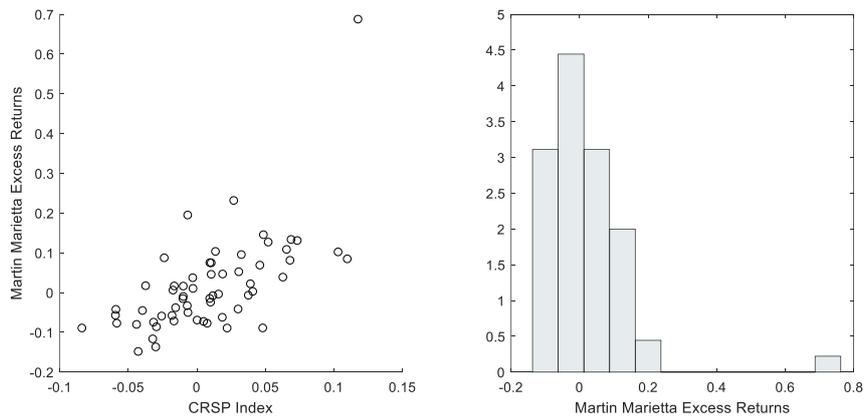

**Figure 2.** (a) Scatter plot of the data set. (b) Histogram of the Martin Marietta excess returns.

In this article, we analyze this data set to illustrate the applicability of the joint location, scale and skewness models of SLN distribution over the joint location, scale and skewness models of the STN distribution. For the comparison of the models, we use the values of the Akaike information criterion (AIC) (Akaike (1973)), the Bayesian information criterion (BIC) (Schwarz (1978)), and the efficient determination criterion (EDC) (Bai et al. (1989))). These criteria have the following form

$$-2\ell(\hat{\boldsymbol{\theta}}) + mc_n,$$

where $\ell(\cdot)$ represents the maximized log-likelihood, $m$ is the number of free parameters to be estimated in the model and $c_n$ is the penalty term. Here, we take $c_n = 2$ for AIC, $c_n = \log(n)$ for BIC and $c_n = 0.2\sqrt{n}$ for EDC.

We give the estimation results in Table 4 for all models. This table contains the estimates, bootstrap standard errors (BSEs) (Efron and Tibshirani (1993)) of estimates based on 500 random samples, the log-likelihood, and the values of AIC, BIC, and EDC. Note that we take the heteroscedastic t model results given in Taylor and Verbyla (2004) as initial values for the parameters of location and scale models. Also, we set $\alpha_0 = \alpha_1 = 0$ as initial values for the parameters of skewness model and the degrees of freedom parameter 3.75. In Figure 3, we show the scatter plot of the data set with the fitted regression lines obtained from the joint location, scale and skewness models of the SLN distribution and the joint location, scale and skewness models of the STN distribution. We observe that the joint location, scale and skewness models of the SLN distribution has better fit than the location, scale and skewness models of the STN distribution according to the information criteria and also Figure 3.

**Table 4.** Estimation results for Martin Marietta data set.

|  |  | Skew t Normal | | Skew Laplace Normal | |
| --- | --- | --- | --- | --- | --- |
|  |  | Estimate | BSE | Estimate | BSE |
| Location model | $\beta_0$ | -0.0349 | 0.0289 | -0.0267 | 0.0227 |
|  | $\beta_1$ | 0.4888 | 1.1387 | 0.7344 | 0.3683 |
| Scale model | $\gamma_0$ | -5.7905 | 0.4418 | -5.8282 | 0.3234 |
|  | $\gamma_1$ | 25.1552 | 15.8670 | 17.2765 | 11.8023 |
| Skewness model | $\alpha_0$ | 0.5614 | 1.6808 | 0.4040 | 1.0923 |
|  | $\alpha_1$ | 19.3303 | 33.9531 | 13.3093 | 9.0474 |
| Information Criteria | $\ell(\hat{\boldsymbol{\theta}})$ | 75.9986 | | **76.9986** | |
|  | AIC | -139.9872 | | **-139.9971** | |
|  | BIC | -125.3268 | | **-127.4311** | |
|  | EDC | -118.3268 | | **-121.4311** | |

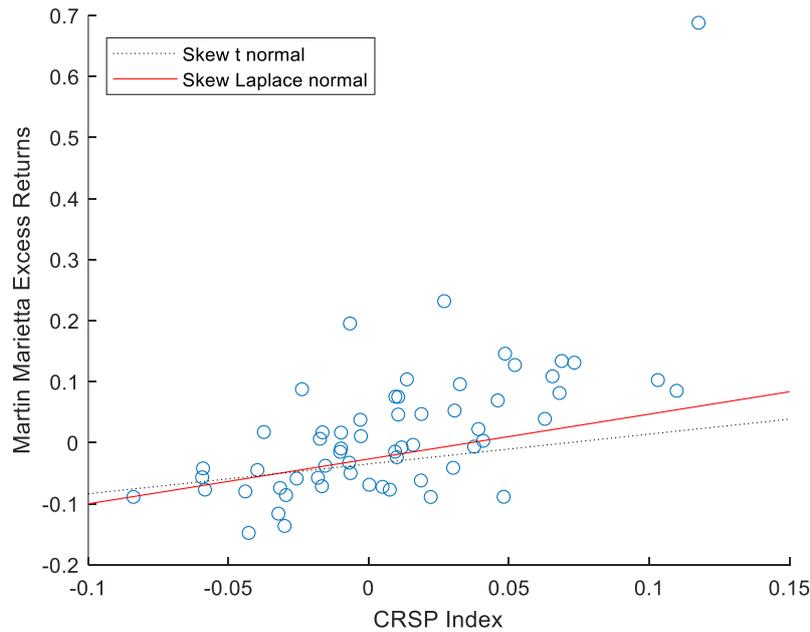

**Figure 3.** The scatterplot of the data set with the fitted regression lines obtained from joint location, scale and skewness models of the SN and SLN distributions.

## 7. Conclusions

In this study, we have proposed the joint location, scale and skewness models of the SLN distribution as an alternative to the joint location, scale and skewness models of the STN distribution. We have obtained the ML estimates via the EM algorithm. We have provided a simulation study to show the estimation performance of the proposed model. We have observed from simulation results that the parameters can be accurately estimated. We have given a real data application to test the applicability of the proposed model and also to compare with the joint, location, scale and skewness models of the STN distribution. We have seen from real data example results that the joint location, scale and skewness models of the SLN distribution gives better fit than the joint, location, scale and skewness models of the STN distribution. Thus, we have concluded that the proposed model can be used as an alternative to the the joint, location, scale and skewness models of the STN distribution for modelling the data sets which have asymmetric and heavy-tailed outcomes.

## References

Akaike, H. 1973. Information theory and an extension of the maximum likelihood principle. Proceeding of the Second International Symposium on Information Theory, B.N. Petrov and F. Caski, eds., 267-281, Akademiai Kiado, Budapest.

Aitkin, M. 1987. Modelling variance heterogeneity in normal regression using GLIM. Applied statistics, 332-339.

Azzalini, A. 1985. A class of distributions which includes the normal ones. Scandinavian Journal of Statistics, 12(2), 171-178.

Azzalini, A. 1986. Further results on a class of distributions which includes the normal ones. Statistica, 46(2), 199-208.

Azzalini, A., Capitanio, A. 2003. Distributions generated by perturbation of symmetry with emphasis on a multivariate skew t distribution. Journal of the Royal Statistical Society, Series B (Statistical Methodology), 65(2), 367-389.

Bai, Z.D., Krishnaiah, P R., Zhao, L.C. 1989. On rates of convergence of efficient detection criteria in signal processing with white noise. IEEE Transactions on Information Theory, 35, 380-388.

DiCiccio, T.J., Monti, A.C. 2004. Inferential aspects of the skew exponential power distribution. Journal of the American Statistical Association, 99(466), 439-450.

Efron, B., Tibshirani, R.J. 1993. An Introduction to the Bootstrap. NewYork: Chapman&Hall.

Engel, J., Huele, A.F. 1996. A generalized linear modeling approach to robust design. Technometrics, 38(4), 365-373.

Gómez, H.W., Venegas, O., Bolfarine, H. 2007. Skew-symmetric distributions generated by the distribution function of the normal distribution. Environmetrics, 18, 395-407.

Harvey, A.C. 1976. Estimating regression models with multiplicative heteroscedasticity. Econometrica: Journal of the Econometric Society, 44(3), 461-465.

Li, H.Q., Wu, L.C. 2014. Joint modelling of location and scale parameters of the skew-normal distribution. Applied Mathematics-A Journal of Chinese Universities, 29(3), 265-272.

Li, H., Wu, L., Ma, T. 2017. Variable selection in joint location, scale and skewness models of the skew-normal distribution, Journal of Systems Science and Complexity, 30(3), 694-709.

Lin, T.I., Wang, Y.J. 2009. A robust approach to joint modeling of mean and scale covariance for longitudinal data. Journal of Statistical Planning and Inference, 139(9), 3013-3026.

Lin, T.I., Wang, W.L. 2011. Bayesian inference in joint modelling of location and scale parameters of the t distribution for longitudinal data. Journal of Statistical Planning and Inference, 141(4), 1543-1553.

Park, R.E. 1966. Estimation with heteroscedastic error terms, Econometrica, 1966, 34, 888.

Schwarz, G. 1978. Estimating the dimension of a model. The Annals of Statistics, 6(2), 461-464.

Smyth, G.K., Verbyla, A.P. 1999. Adjusted likelihood methods for modelling dispersion in generalized linear models. Environmetrics, 10(6), 695-709.

Taylor, J., Verbyla, A. 2004. Joint modelling of location and scale parameters of the t distribution. Statistical Modelling, 4(2), 91-112.

Verbyla, A.P. 1993. Modelling variance heterogeneity: residual maximum likelihood and diagnostics. Journal of the Royal Statistical Society. Series B (Methodological), 55(2), 493-508.

Wu, L., Li, H. 2012. Variable selection for joint mean and dispersion models of the inverse Gaussian distribution. Metrika, 75(6), 795-808.

Wu, L., Tian, G. L., Zhang, Y. Q., Ma, T. 2017. Variable selection in joint location, scale and skewness models with a skew-t-normal distribution. Statistics and Its Interface, 10(2), 217-227.

Zhao, W., Zhang, R. 2015. Variable selection of varying dispersion student-t regression models. Journal of Systems Science and Complexity, 28(4), 961-977.

**Appendix**

Using the objective function given in (21), we obtain the score function

$$G(\boldsymbol{\theta}) = \frac{\partial Q(\boldsymbol{\theta}; \hat{\boldsymbol{\theta}})}{\partial \boldsymbol{\theta}} = \left(G_1^T(\boldsymbol{\beta}), G_2^T(\boldsymbol{\gamma}), G_3^T(\boldsymbol{\alpha})\right)^T,$$

where

$$G_1(\boldsymbol{\beta}) = \sum_{i=1}^{n} \frac{(y_i - \boldsymbol{x}_i^T \boldsymbol{\beta}) \boldsymbol{x}_i}{e^{\boldsymbol{z}_i^T \boldsymbol{\gamma}}} \left(\hat{v}_i + (\boldsymbol{w}_i^T \boldsymbol{\alpha})^2\right) - \sum_{i=1}^{n} \frac{(\boldsymbol{w}_i^T \boldsymbol{\alpha}) \boldsymbol{x}_i \hat{u}_{1i}}{e^{\boldsymbol{z}_i^T \boldsymbol{\gamma}/2}},$$

$$G_2(\boldsymbol{\gamma}) = -\frac{1}{2}\sum_{i=1}^{n} \boldsymbol{z}_i + \frac{1}{2}\sum_{i=1}^{n} \frac{(y_i - \boldsymbol{x}_i^T \boldsymbol{\beta})^2 \boldsymbol{z}_i}{e^{\boldsymbol{z}_i^T \boldsymbol{\gamma}}} \left(\hat{v}_i + (\boldsymbol{w}_i^T \boldsymbol{\alpha})^2\right) - \frac{1}{2}\sum_{i=1}^{n} \frac{(\boldsymbol{w}_i^T \boldsymbol{\alpha})(y_i - \boldsymbol{x}_i^T \boldsymbol{\beta}) \boldsymbol{z}_i \hat{u}_{1i}}{e^{\boldsymbol{z}_i^T \boldsymbol{\gamma}/2}},$$

$$G_3(\boldsymbol{\alpha}) = \sum_{i=1}^{n} \frac{(y_i - \boldsymbol{x}_i^T\boldsymbol{\beta})\boldsymbol{w}_i \hat{u}_{1i}}{e^{\boldsymbol{z}_i^T\boldsymbol{\gamma}/2}} - \sum_{i=1}^{n} \frac{(\boldsymbol{w}_i^T\boldsymbol{\alpha})(y_i - \boldsymbol{x}_i^T\boldsymbol{\beta})\boldsymbol{w}_i}{e^{\boldsymbol{z}_i^T\boldsymbol{\gamma}}},$$

and observed Fisher information matrix

$$H(\boldsymbol{\theta}) = \frac{\partial^2 Q(\boldsymbol{\theta};\widehat{\boldsymbol{\theta}})}{\partial \boldsymbol{\theta} \partial \boldsymbol{\theta}^T} = \begin{bmatrix} \frac{\partial^2 Q(\boldsymbol{\theta};\widehat{\boldsymbol{\theta}})}{\partial \boldsymbol{\beta} \partial \boldsymbol{\beta}^T} & \frac{\partial^2 Q(\boldsymbol{\theta};\widehat{\boldsymbol{\theta}})}{\partial \boldsymbol{\beta} \partial \boldsymbol{\gamma}^T} & \frac{\partial^2 Q(\boldsymbol{\theta};\widehat{\boldsymbol{\theta}})}{\partial \boldsymbol{\beta} \partial \boldsymbol{\alpha}^T} \\ \frac{\partial^2 Q(\boldsymbol{\theta};\widehat{\boldsymbol{\theta}})}{\partial \boldsymbol{\gamma} \partial \boldsymbol{\beta}^T} & \frac{\partial^2 Q(\boldsymbol{\theta};\widehat{\boldsymbol{\theta}})}{\partial \boldsymbol{\gamma} \partial \boldsymbol{\gamma}^T} & \frac{\partial^2 Q(\boldsymbol{\theta};\widehat{\boldsymbol{\theta}})}{\partial \boldsymbol{\gamma} \partial \boldsymbol{\alpha}^T} \\ \frac{\partial^2 Q(\boldsymbol{\theta};\widehat{\boldsymbol{\theta}})}{\partial \boldsymbol{\alpha} \partial \boldsymbol{\beta}^T} & \frac{\partial^2 Q(\boldsymbol{\theta};\widehat{\boldsymbol{\theta}})}{\partial \boldsymbol{\alpha} \partial \boldsymbol{\gamma}^T} & \frac{\partial^2 Q(\boldsymbol{\theta};\widehat{\boldsymbol{\theta}})}{\partial \boldsymbol{\alpha} \partial \boldsymbol{\alpha}^T} \end{bmatrix},$$

where

$$\frac{\partial^2 Q(\boldsymbol{\theta};\widehat{\boldsymbol{\theta}})}{\partial \boldsymbol{\beta} \partial \boldsymbol{\beta}^T} = -\sum_{i=1}^{n} \frac{\boldsymbol{x}_i \boldsymbol{x}_i^T}{e^{\boldsymbol{z}_i^T\boldsymbol{\gamma}}} \hat{v}_i - \sum_{i=1}^{n} \frac{(\boldsymbol{w}_i^T\boldsymbol{\alpha})^2}{e^{\boldsymbol{z}_i^T\boldsymbol{\gamma}}} \boldsymbol{x}_i \boldsymbol{x}_i^T,$$

$$\frac{\partial^2 Q(\boldsymbol{\theta};\widehat{\boldsymbol{\theta}})}{\partial \boldsymbol{\beta} \partial \boldsymbol{\gamma}^T} = -\sum_{i=1}^{n} \frac{(y_i - \boldsymbol{x}_i^T\boldsymbol{\beta})\boldsymbol{x}_i \boldsymbol{z}_i^T}{e^{\boldsymbol{z}_i^T\boldsymbol{\gamma}}} \left(\hat{v}_i + (\boldsymbol{w}_i^T\boldsymbol{\alpha})^2\right) + \frac{1}{2}\sum_{i=1}^{n} \frac{(\boldsymbol{w}_i^T\boldsymbol{\alpha})\boldsymbol{x}_i \boldsymbol{z}_i^T \hat{u}_{1i}}{e^{\boldsymbol{z}_i^T\boldsymbol{\gamma}/2}},$$

$$\frac{\partial^2 Q(\boldsymbol{\theta};\widehat{\boldsymbol{\theta}})}{\partial \boldsymbol{\beta} \partial \boldsymbol{\alpha}^T} = -\sum_{i=1}^{n} \frac{\boldsymbol{x}_i \boldsymbol{w}_i^T \hat{u}_{1i}}{e^{\boldsymbol{z}_i^T\boldsymbol{\gamma}/2}} + 2\sum_{i=1}^{n} \frac{(\boldsymbol{w}_i^T\boldsymbol{\alpha})(y_i - \boldsymbol{x}_i^T\boldsymbol{\beta})\boldsymbol{x}_i \boldsymbol{w}_i^T}{e^{\boldsymbol{z}_i^T\boldsymbol{\gamma}}},$$

$$\frac{\partial^2 Q(\boldsymbol{\theta};\widehat{\boldsymbol{\theta}})}{\partial \boldsymbol{\gamma} \partial \boldsymbol{\beta}^T} = -\sum_{i=1}^{n} \frac{(y_i - \boldsymbol{x}_i^T\boldsymbol{\beta})\boldsymbol{z}_i \boldsymbol{x}_i^T}{e^{\boldsymbol{z}_i^T\boldsymbol{\gamma}}} \left(\hat{v}_i + (\boldsymbol{w}_i^T\boldsymbol{\alpha})^2\right) + \frac{1}{2}\sum_{i=1}^{n} \frac{(\boldsymbol{w}_i^T\boldsymbol{\alpha})\boldsymbol{z}_i \boldsymbol{x}_i^T \hat{u}_{1i}}{e^{\boldsymbol{z}_i^T\boldsymbol{\gamma}/2}},$$

$$\frac{\partial^2 Q(\boldsymbol{\theta};\widehat{\boldsymbol{\theta}})}{\partial \boldsymbol{\gamma} \partial \boldsymbol{\gamma}^T} = -\frac{1}{2}\sum_{i=1}^{n} \frac{(y_i - \boldsymbol{x}_i^T\boldsymbol{\beta})^2 \boldsymbol{z}_i \boldsymbol{z}_i^T}{e^{\boldsymbol{z}_i^T\boldsymbol{\gamma}}} \left(\hat{v}_i + (\boldsymbol{w}_i^T\boldsymbol{\alpha})^2\right) + \frac{1}{4}\sum_{i=1}^{n} \frac{(\boldsymbol{w}_i^T\boldsymbol{\alpha})(y_i - \boldsymbol{x}_i^T\boldsymbol{\beta})\boldsymbol{z}_i \boldsymbol{z}_i^T \hat{u}_{1i}}{e^{\boldsymbol{z}_i^T\boldsymbol{\gamma}/2}},$$

$$\frac{\partial^2 Q(\boldsymbol{\theta};\widehat{\boldsymbol{\theta}})}{\partial \boldsymbol{\gamma} \partial \boldsymbol{\alpha}^T} = -\frac{1}{2}\sum_{i=1}^{n} \frac{(y_i - \boldsymbol{x}_i^T\boldsymbol{\beta})\boldsymbol{z}_i \boldsymbol{w}_i^T \hat{u}_{1i}}{e^{\boldsymbol{z}_i^T\boldsymbol{\gamma}/2}} + \sum_{i=1}^{n} \frac{(\boldsymbol{w}_i^T\boldsymbol{\alpha})(y_i - \boldsymbol{x}_i^T\boldsymbol{\beta})^2 \boldsymbol{z}_i \boldsymbol{w}_i^T}{e^{\boldsymbol{z}_i^T\boldsymbol{\gamma}}},$$

$$\frac{\partial^2 Q(\boldsymbol{\theta};\widehat{\boldsymbol{\theta}})}{\partial \boldsymbol{\alpha} \partial \boldsymbol{\beta}^T} = -\sum_{i=1}^{n} \frac{\boldsymbol{w}_i \boldsymbol{x}_i^T \hat{u}_{1i}}{e^{\boldsymbol{z}_i^T\boldsymbol{\gamma}/2}} + 2\sum_{i=1}^{n} \frac{(\boldsymbol{w}_i^T\boldsymbol{\alpha})(y_i - \boldsymbol{x}_i^T\boldsymbol{\beta})\boldsymbol{w}_i \boldsymbol{x}_i^T}{e^{\boldsymbol{z}_i^T\boldsymbol{\gamma}}},$$

$$\frac{\partial^2 Q(\boldsymbol{\theta};\widehat{\boldsymbol{\theta}})}{\partial \boldsymbol{\alpha} \partial \boldsymbol{\gamma}^T} = -\frac{1}{2}\sum_{i=1}^{n} \frac{(y_i - \boldsymbol{x}_i^T\boldsymbol{\beta})\boldsymbol{w}_i \boldsymbol{z}_i^T \hat{u}_{1i}}{e^{\boldsymbol{z}_i^T\boldsymbol{\gamma}/2}} + \sum_{i=1}^{n} \frac{(\boldsymbol{w}_i^T\boldsymbol{\alpha})(y_i - \boldsymbol{x}_i^T\boldsymbol{\beta})^2 \boldsymbol{w}_i \boldsymbol{z}_i^T}{e^{\boldsymbol{z}_i^T\boldsymbol{\gamma}}},$$

$$\frac{\partial^2 Q(\boldsymbol{\theta};\widehat{\boldsymbol{\theta}})}{\partial \boldsymbol{\alpha} \partial \boldsymbol{\alpha}^T} = -\sum_{i=1}^{n} \frac{(y_i - \boldsymbol{x}_i^T\boldsymbol{\beta})^2 \boldsymbol{w}_i \boldsymbol{w}_i^T}{e^{\boldsymbol{z}_i^T\boldsymbol{\gamma}}}.$$